\documentclass[12pt]{article}
\usepackage{amsfonts} \usepackage{amssymb} \usepackage{textcomp}

\renewcommand{\phi}{\varphi} \newcommand{\be}{\begin{equation}}
\newcommand{\ee}{\end{equation}}
\newcommand{\ba}{\begin{eqnarray}}
\newcommand{\ea}{\end{eqnarray}}
\newcommand{\ban}{\begin{eqnarray*}}
\newcommand{\ean}{\end{eqnarray*}}

\newcommand{\nul}{{\bf0}} \newcommand{\cupm}{\mathop{\cup}}
 
 \newcommand{\rd}{{\Bbb R}^d}
  
\newcommand{\zd}{{\Bbb Z}^{d}} 
 \renewcommand{\r}{{\Bbb R}}
\newcommand{\z} {{\Bbb Z}} \newcommand{\n} {{\Bbb N}}
 
 \newcommand{\ddd}{,\dots,}
\renewcommand{\lll}{\left(} \newcommand{\rrr}{\right)}
\newcommand{\ex}[1]{e^{2\pi i{#1}}} \newcommand{\exm}[1]{e^{-2\pi
i{#1}}} \newcommand{\sml}[3]{\sum\limits_{{#1}={#2}}^{#3}}

\newtheorem{theo}{Theorem} \newtheorem{lem}[theo]{Lemma}
\newtheorem {prop} [theo] {Proposition} \newtheorem {coro} [theo]
{Corollary}

\begin {document} \centerline{\Large\bf{Decompositions of
Trigonometric} }

\centerline{\Large\bf { Polynomials with  Applications to }}

\centerline{\Large\bf { Multivariate Subdivision Schemes
\footnote {
This research was supported by The Hermann
Minkowski Center for Geometry at Tel-Aviv University.
The
second author is also supported  by DFG Project 436 RUS 113/951 and
Grant 09-01-00162 of RFBR.}
}}
\bigskip \centerline{\large Nira Dyn
and Maria Skopina } \bigskip

\centerline{Tel-Aviv University, e-mail: niradyn@post.tau.ac.il }
\centerline{and }
\centerline{St.Petersbutg State University, e-mail: skopina@MS1167.spb.edu}

\begin{abstract}
We study  multivariate trigonometric polynomials,
satisfying a set of constraints
close  to the known Strung-Fix conditions.
Based on the polyphase representation of these
polynomials relative to a general dilation matrix, we develop a simple
constructive method  for a special type of decomposition of such  polynomials.
These decompositions are of interest to the analysis of
convergence and smoothness of multivariate subdivision schemes
associated with general dilation matrices.
We apply these decompositions, by verifying sufficient conditions for the
convergence and smoothness of multivariate scalar subdivision schemes,
proved here. For the convergence analysis
our sufficient conditions apply to arbitrary
dilation matrices, while the previously known necessary and sufficient
conditions are relevant only in case  of
dilation matrices with a self similar tiling.
For the analysis of smoothness, we state and prove two theorems on
multivariate  matrix subdivision schemes, which lead to sufficient conditions for
$C^1$ limits of scalar multivariate subdivision schemes associated with
isotropic dilation matrices. Although similar results are stated in the
literature, we give here detailed proofs of the results, which we could
not find elsewhere.
\end{abstract}

 MSC  41A30, 42C40

\subsection{1. Introduction}

In this  paper we study  multivariate trigonometric polynomials
(masks), satisfying a set of constraints, which
is close  to the known Strung-Fix conditions. We are interested in decomposing such a
polynomial multiplied by a factor of the form $1-\ex{(x,{\bf
e}_j)}$,  into a sum of $d$ trigonometric polynomials, satisfying
the same type of constraints, but of one order less, each
multiplied by a different factor of the form $1-\ex{(M^*x,{\bf
e}_k)}$.

This type of decompositions is of interest in the analysis of
convergence and smoothness of multivariate subdivision schemes.
The algebraic approach to the analysis of such schemes is
investigated in~\cite{0} relative to the dilation matrix $2I$, and
in~\cite{10}, \cite{13}, \cite{14}, relative to general dilation
matrices.

The analogous decompositions for multivariate algebraic
polynomials are based on the rich structure of ideals of
multivariate polynomials. The trigonometric decompositions are
based on the rather simple idea of the representation of a
trigonometric polynomial in terms of its unique polyphase
trigonometric polynomials relative to  a dilation matrix~\cite{2}.

The outline of the paper is as follows: in Section 2 we introduce
our notation and the set of constraints. We also bring basic
results related to polyphase representations of trigonometric
polynomials. Properties of trigonometric polynomials satisfying
the set of constraints of order zero relative to a general
dilation matrix are derived in Section 3, in terms of their
polyphase trigonometric polynomials.

These results are used in Section 4 to prove the decomposition of
trigonometric polynomials, satisfying the set of constraints of
order zero. An  algorithm for the computation of this
decomposition is presented. The decomposition for order zero is
used to derive that for order one. The general case of order
greater than one, which is much more involved, is also
investigated, and an algorithm for this case is given.

Section 5 gives a characterization of a trigonometric polynomial
satisfying a set of constraints of a general order, in terms of
the values at the origin of the derivatives of its polyphase
trigonometric polynomials. This result provides a tool for the
construction of such polynomials.

Applications of the decompositions of Section 4 to the convergence
and smoothness analysis of multivariate subdivision schemes
associated with general dilation matrices are presented in
Section~6.
The application of the decompositions to the convergence analysis is
based on a new result giving sufficient conditions for the convergence
of a scalar multivariate subdivision scheme associated with an arbitrary
dilation matrix.
The  result in the literature giving necessary and sufficient conditions
for convergence, is limited to dilation matrices
having a self similar tiling~\cite{10}.
For the analysis of smoothness, we state and prove two theorems on
multivariate matrix subdivision schemes. First we give a detailed proof of a
sufficient condition for $C^1$  limit functions of such a scheme,
associated with an isotropic dilation matrix. Although  this
condition is stated in ~\cite{14}, the sketch of a proof given
there, is valid only for dilation matrices which are multiples of
the identity matrix. Then we state and prove a condition, in terms
of the decompositions, that allows us to check the sufficient
condition for $C^1$.

\subsection {2. Notation and  preliminary information}

Let $ \n $ be the set of positive integers, $ \rd $ denotes the $d
$-dimensional Euclidean space, $x = (x_1\ddd x_d) $, $y = (y_1\ddd
y_d) $ are its elements (vectors), $ (x, y) = x_1y_1
+\dots+x_dy_d$,
$|x|=\max\limits_{j=1\ddd d}|x_j|$, ${\bf e}_j=(0\ddd1\ddd0)$ is
the $j$-th unit vector in $\rd$, \ $ \nul = (0\ddd0) \in \rd $; $
\zd $ is the integer lattice in $ \rd $. For $x,y\in\rd$, we
write $x>y$ if $x_j>y_j$, $j=1\ddd d$; $ \zd_+ = \{x\in Z^d:\
x\ge\nul \} $. If $\alpha, \beta\in\zd_+$, $a,b\in\rd$, we set
$\alpha!=\prod\limits_{j=1}^d\alpha_j!$,
$\lll\alpha\atop\beta\rrr= \frac{\alpha!}{\beta!(\alpha-\beta)!}$,
$a^b=\prod\limits_{j=1}^d{a_j}^{b_j}$,
 $[\alpha]=\sml j1d \alpha_j$,
$D^\alpha f=\frac{\partial^{[\alpha]}f}
{\partial^{\alpha_1}x_1\dots\partial^{\alpha_d}x_d}$;
$\delta_{ab}$ denotes the Kronecker delta; $z:=(z_1\ddd z_d)$,
where $z_k:=\ex{x_k}$;
$\ell_\infty^N(\zd)=\ell_\infty^N:=\{f=\{f_\alpha\}_{\alpha\in\zd}:\
\
 f_\alpha\in\r^N, |f_\alpha|\le C\}$, $\|f\|_{\ell_\infty^N}=\|f\|_\infty:=
 \sup_{\alpha\in\zd}|f_\alpha|$,
 $\ell_\infty:=\ell_\infty^1$.

Let $ M $ be a non-degenerate $ d \times d $ integer matrix
 whose eigenvalues are bigger than 1 in module,
$M^* $ is the transpose of $M$, $I_d$  and ${\Bbb O}_d$ denote
respectively the unit and the zero $ d \times d$ matrices. We  say
that the numbers $k, n\in \zd $ are congruent modulo~ $M $ (write
$k\equiv n  \pmod {M} $) if $k-n=M\ell $, $ \ell\in\zd $. The
integer lattice $ \zd $ is splitted into cosets with respect to
the introduced relation of congruence. The number of cosets is
equal to
 $ | \det M | $ (see, e.g., \cite [p. 107] {73}).
Let us take  an arbitrary representative from each coset, call
them digits and denote the set of digits by $D (M) $. Throughout
the paper we  consider that such a matrix $M$ ({\em dilation
matrix}) is fixed, $m=|\det M|$, $D (M) =\{s_0\ddd s_{m-1}\}$, $D
(M^*) =\{s^*_0\ddd s^*_{m-1}\}$, $s_0=s^*_0=\nul$,
$r_k=M^{-1}s_k$,
 $k=0\ddd m-1$.

If  $A$ is a $N\times N'$ matrix with entries $a_{jk}$, then
$A^{[n]}$ denotes its $n$-th Kronecker power, i.e., $$
A^{[0]}:=1,\ \ \ A^{[1]}:=A,\ \ \ A^{[n+1]}:= \lll\matrix{
a_{11}A^{[n]}&\dots&a_{1N'}A^{[n]} \cr \vdots&\ddots&\vdots \cr \cr
a_{N1}A^{[n]}&\dots&a_{NN'}A^{[n]}} \rrr. $$

\noindent {\bf Proposition  A} {\em The matrix
$\left\{\frac1{\sqrt m}\,\ex{(r_k, s_l^*)}\right\}_{k,l=0}^{m-1}$
is unitary. In particular,} \be \frac1{\sqrt
m}\sum\limits_{k=0}^{m-1}\ex{(r_k, s_l^*)}=\delta_{0l}.
\ee

A proof of this statement can be found in~\cite{NPS} or
in~\cite{1}.

We will consider $1$-periodic trigonometric polynomials in $d$
variables $$ t(x)=\sum\limits_{n\in\zd}\widehat t(n)\ex{(n,x)}, $$
where the set $\{n\in\zd : \widehat t(n)\ne 0\}$ is finite.

For any  trigonometric polynomial $t$, there exists a unique set
of   trigonometric polynomials $\tau_{\nu}$, $\nu=0\ddd m-1$,
(polyphase functions of $t$) (see, e.g.,~\cite{NPS}) such that \be
t(x)=\sml\nu0{m-1}\ex{(s_\nu,x)}\tau_{\nu}(M^*x). \label{1} \ee It
is clear that $$ \tau_{\nu}(x)=\sum\limits_{n\in\zd}\widehat
t(Mn+s_\nu)\ex{(n,x)}. $$

For any $n\in\z_+$, denote by ${\cal Z}^n$ the set of
trigonometric polynomials $t$ satisfying the following
''zero-condition''
 of order~$n$:
$D^\beta t({M^*}^{-1}x)|_{x=s}=0$ for all $s\in D(M^*)$,
$s\ne\nul$, and for all $\beta\in\zd_+$, $[\beta]\le n$. It is
well known that this ''zero-condition'' is equivalent to the
so-called ''sum-rule'' and close to the Strang-Fix condition of
the same order. It will be convenient for us to define ${\cal
Z}^{-1}$ as the set of all trigonometric polynomials.

\subsection {3. Auxiliary results}

\begin{prop} A trigonometric polynomial $t$ belongs to ${\cal
Z}^0$ if and only if \be \tau_\nu(\nul)=\frac{t(\nul)}{m},\ \ \
\nu=0\ddd m-1, \label{01} \ee where $\tau_{0}\ddd\tau_{m-1}$ are
the  polyphase functions of $t$. \label{p1} \end{prop}

{\bf Proof.} Let $s\in D(M^*)$, using~(\ref{1}) we have $$
t({M^*}^{-1}s)=\sml \nu 0{m-1}\ex{(r_\nu,s)}\tau_{\nu}(s)= \sml
\nu 0{m-1}\ex{(r_\nu,s)}\tau_{\nu}(\nul). $$ So, the relation
$t\in {\cal Z}^0$ can be rewritten as $$ \sml \nu
0{m-1}\ex{(r_\nu,s_l^*)}\tau_{\nu}(\nul)=t(\nul)\delta_{\nul l}, \
\ \  l=0\ddd m-1. $$ Consider these equalities as a linear system
with unknowns $\tau_0(\nul)\ddd\tau_{m-1}(\nul)$. Due to
Proposition A, the system has a unique solution given
by~(\ref{01}). $\Diamond$

\begin{lem} Let $t$ be a trigonometric polynomial, $t(\nul)=0$,
then $$ t(x)=\sml k1d t_k(x)(1-\ex{(x,{\bf e}_k)}), $$ where
$t_k$, $k=1\ddd d$, are trigonometric polynomials. \label{l1}
\end{lem}

{\bf Proof.} There exists $N\in\zd$ such that
$t(x)\ex{(N,x)}=p(z)$, where  $z=\exp(2\pi ix)$ and $p$ is an
algebraic polynomial. It is clear that $P(1\ddd1)=0$. By Taylor
formula it follows that $$ p(z)=\sml k1d p_k(z)(1-z_k). $$ where
$p_k$, $k=1\ddd d$, are algebraic polynomials. It remains to set
$t_k(x):=p_k(z)\exm{(N,x)}. \Diamond$

\begin{lem} Let $t, \tilde t\in {\cal Z}^0$, $t(\nul)=\tilde
t(\nul)$, then
\be
t(x)-\tilde t(x)=\sml k1d t_k(x)\lll
1-\ex{(M^*x,{\bf e}_k)}\rrr,
\label{02}
\ee
 where $t_k$, $k=1\ddd
d$, are trigonometric polynomials. \label{l2} \end{lem}

{\bf Proof.} Using~(\ref{1}) for $t, \tilde t$, we have \be t
(x)-\tilde t(x)=\sml \nu
0{m-1}\ex{(s_\nu,x)}(\tau_{\nu}(M^*x)-\tilde \tau_{\nu}(M^*x)).
\label{2} \ee Due to Proposition~\ref{p1},
$\tau_{\nu}(\nul)-\tilde \tau_{\nu}(\nul)=0$, $\nu=1\ddd m-1$.
Hence, by Lemma~\ref{l1}, \be \tau_{\nu}(y)-\tilde
\tau_{\nu}(y)=\sml l1d \tau_{\nu k}(y)\lll 1-\ex{(y,{\bf
e}_k)}\rrr, \label{3} \ee where $\tau_{\nu k}$, $k=1\ddd d$, are
trigonometric polynomials. It remains to set $y=M^*x$ and
combine~(\ref{3}) with~(\ref{2}). $\Diamond$

\begin{coro} Let $t\in {\cal Z}^0$, $r\in\zd$, $\tilde
t:=\ex{(r,\cdot)}t$, then~(\ref{02}) holds. \label{l3} \end{coro}

\begin{lem} Let $t\in {\cal Z}^{-1}$, $r\in\zd$, $\tilde
t:=\ex{(r,\cdot)}t$, then the $k$-th polyphase function of $\tilde
t$ is given by $\tilde \tau_{k}(x):=\ex{(n_\nu,x)}\tau_\nu(x)$,
where $\nu$ is the number of the unique digit $s_\nu\in D(M)$ such
that $s_\nu+r=s_{k}+Mn_\nu$, $n_\nu\in\zd$. \label{l5} \end{lem}

{\bf Proof.} By~(\ref{1}), $$ \tilde t(x)=\ex{(r,x)}\sml
\nu0{m-1}\ex{(s_\nu,x)}\tau_{\nu}(M^*x)= \sml
\nu0{m-1}\ex{(s_\nu+r,x)}\tau_{\nu}(M^*x). $$ It is clear that
$s_\nu+r=s_{l_\nu}+Mn_\nu$, $n_\nu\in\zd$, $\{l_0\ddd
l_{m-1}\}=\{0\ddd m-1\}$. So, $$ \tilde t(x)=\sml
\nu0{m-1}\ex{(s_{l_\nu},x)}\ex{(n_\nu,M^*x)}\tau_{\nu}(M^*x), $$
and the trigonometric polynomial $\tilde
\tau_{l_\nu}(x):=\ex{(n_\nu,x)}\tau_\nu(x)$ is the $l_\nu$-th
polyphase function of $\tilde t$. $\Diamond$

\begin{lem} Let $n\in\zd$, $ 1-\ex{(M^{-1}n,{\bf e}_j)}=0$ for
each $j=1\ddd d$. Then $n\equiv\nul\pmod M$. \label{l4} \end{lem}

{\bf Proof.} Since $\ex n=1$ if and only if $n\in\z$, we have
$(M^{-1}n)_j=(M^{-1}n,{\bf e}_j)\in\z$ for each $j=1\ddd d$.
Hence, $M^{-1}n=l\in\zd$, i.e. $n=Ml$.~$\Diamond$

\subsection {4. Decomposition of masks}

In this section we prove the existence of a decomposition of a
trigonometric polynomial $t\in {\cal Z}^n$ of the form
$$(1-\ex{(x,{\bf e}_k)})t(x)=\sml j1d t_{jk}(x)(1-\ex{(M^*x,{\bf
e}_j)}),$$ with $t_{jk}\in {\cal Z}^{n-1}\ ,j=1\ddd d$, for all
$k=1\ddd d$, and present algorithms for its computation. We start
with the simpler case $n=0,1$, and then do the general case $n>1$,
which is much more complicated.

\subsubsection {4.1 The case $t\in {\cal  Z}^n,\ n=0,1$}

\begin{prop}
 Let $t\in {\cal Z}^0$, then \be (1-\ex{(x,{\bf
e}_k)})t(x)=\sml j1d t_{jk}(x)(1-\ex{(M^*x,{\bf e}_j)}),
\label{001} \ee where $t_{jk}$, $j=1\ddd d$, are trigonometric
polynomials. Conversely, if~(\ref{001}) holds for a trigonometric
polynomial $t$, then $t\in {\cal  Z}^0$ and
$t_{jk}(\nul)=(M^{-1})_{jk}\,t(\nul)$.
\label{p2}
\end{prop}

{\bf Proof.} If  $t\in {\cal  Z}^0$, then~(\ref{001})
 follows immediately form Corollary~\ref{l3}. Now let $t$ be an arbitrary trigonometric
 polynomial such that~(\ref{001}) holds,  $s\in D(M^*)$, $s\ne\nul$. It follows from~(\ref{001}) that
$$
 t({M^*}^{-1}s)\lll 1-\ex{({M^*}^{-1}s,{\bf e}_k))}\rrr=
\\
\sml j1d t_{jk}({M^*}^{-1}s)(1-\ex{(s,{\bf e}_j)})=0. $$ Taking
into account Lemma~\ref{l4}, we obtain $t({M^*}^{-1}s)=0$, which
means that $t\in {\cal  Z}^0$.

 Rewrite identity~(\ref{001}) as follows
$$ (1-\ex{(x,M^{-1}{\bf e}_k)})t({M^*}^{-1}x)=\sml j1d
t_{jk}({M^*}^{-1}x)(1-\ex{(x,{\bf e}_j)}), $$ Differentiating by
$x_l$, we have \ban -2\pi i\,(M^{-1}{\bf e}_k,{\bf
e}_l)\ex{(x,M^{-1}{\bf e}_k)}) t({M^*}^{-1}x)+
(1-\ex{(x,M^{-1}{\bf e}_k)})\frac{\partial}{\partial
x_l}t({M^*}^{-1}x)=
\\
\sml j1d \frac{\partial}{\partial
x_l}t_{jk}({M^*}^{-1}x)(1-\ex{(x,{\bf e}_j)})- 2\pi i\,
t_{lk}({M^*}^{-1}x)\,\ex{(x,{\bf e}_l)}. \ean Substituting
$x=\nul$, we obtain $$ -2\pi i\,(M^{-1}{\bf e}_k,{\bf e}_l)
t(\nul)= -2\pi i\,t_{lk}(\nul). $$ It remains to note that $
(M^{-1}{\bf e}_k,{\bf e}_l)=({M}^{-1})_{lk}.\Diamond $

Analyzing the proofs of  Lemmas~\ref{l2} and \ref{l5} it is not
difficult to describe an algorithm for finding the polynomials
$t_{jk}$, $j,k=1\ddd d$, in the decomposition~(\ref{001}). We
assume that the polyphase functions $\tau_\nu$ of $t$ are given,
and the algorithm derives the polyphase functions $\tau_{jk\nu}$
of each polynomial $t_{jk}$.

\vspace{.5cm} \noindent {\bf ALGORITHM 1}

{\bf Input:} \ \ \ \ \  $\{\widehat\tau_\nu(l),  \ l\in\zd, \
\nu=0\ddd m-1\}$.

{\bf Output:} \ \ $\{\widehat\tau_{jk\nu}(l), \  l\in\zd,\
\nu=0\ddd m-1, \ j,k=1\ddd d\}$.

{\bf Step 0.} \ \ For each $\nu=0\ddd m-1$ find the sets
$$\Omega_\nu:=\{l\in\zd:\ \widehat\tau_\nu(l)\ne0\}.$$

{\bf For each }$k=1\ddd d$ {\bf do}

\hspace{.1cm} {\bf For each} $\nu=0\ddd m-1$ {\bf do}

\begin{quote} {\bf Step 1.} Compute $M^{-1}({\bf
e}_k-s_\nu+s_n)=:l_n$ for all $n=0\ddd m-1$ and denote by $n^*$
the unique $n$ such that $l_n\in\zd$. \end{quote}

\begin{quote} {\bf Step 2.} Find the set
$\tilde\Omega_{k\nu}:=\{l=r+l_{n^*}, r\in \Omega_{n^*}\}$, and for
each $l\in \tilde\Omega_{k\nu}$ put
$\widehat{\tilde\tau}_\nu(l):={\widehat\tau}_{n^*}(l-l_{n^*})$, $$
p_{k\nu}(z):=\sum\limits_{l\in\Omega_\nu\cup{\tilde\Omega}_{k\nu}}
({\widehat\tau}_\nu(l)-\widehat{\tilde\tau}_\nu(l))z^l. $$
\end{quote}

\begin{quote} {\bf Step 3.} Put \ban
&&\tau_{1k\nu}(x)=\frac1{1-z_1}\lll p_{k\nu}(z)-p_{k\nu}(1,z_2\ddd
z_d) \rrr,
\\
&&\tau_{2k\nu}(x)=\frac1{1-z_2}\lll p_{k\nu}(1,z_2,z_3\ddd z_d)
-p_{k\nu}(1,1,z_3\ddd z_d)\rrr,
\\
&&\ldots\ldots\ldots\ldots\ldots\ldots\ldots\ldots\ldots\ldots\ldots\ldots\ldots\ldots\ldots\ldots\ldots
\ldots\ldots
\\
&&\tau_{dk\nu}(x)=\frac1{1-z_d}\lll p_{k\nu}(1,1\ddd 1,z_d)
-p_{k\nu}(1,1\ddd 1,1)\rrr.\Diamond \ean \end{quote}

\begin{theo} If  $t\in {\cal  Z}^1$, then  in any
decomposition~(\ref{001}) the trigonometric polynomials $t_{jk}$,
$j,k=1\ddd d$,  belong to ${\cal  Z}^0$. \label{t1i} \end{theo}

{\bf Proof.}
Let $l=1\ddd d$, $s\in D(M^*)$, $s\ne\nul$. Since $t\in {\cal
Z}^{(1)}$, \be \frac{\partial}{\partial x_l}\lll \lll
1-\ex{({M^*}^{-1}x,{\bf e}_k))}\rrr
t({M^*}^{-1}x)\rrr\Bigg|_{x=s}=0, \ \ \ k=1\ddd d. \label{28} \ee
On the other hand, it follows from~(\ref{001}) that for each
$k=1\ddd d$ we have \ban \frac{\partial}{\partial x_l}\lll\lll
1-\ex{({M^*}^{-1}x,{\bf e}_k))}\rrr
t({M^*}^{-1}x)\rrr\Bigg|_{x=s}= \hspace{3cm}
\\
\frac{\partial}{\partial x_l}\lll\sml j1d
t_{jk}({M^*}^{-1}x)(1-\ex{x_j})\rrr\Bigg|_{x=s}=
\\
\sml j1d t_{jk}({M^*}^{-1}s)\frac{\partial}{\partial
x_l}(1-\ex{x_j})\Bigg|_{x=s}= -2\pi i\,t_{lk}({M^*}^{-1}s). \ean
So, we proved that $t_{lk}({M^*}^{-1}s)=0$ for all  $s\in D(M^*)$,
$s\ne\nul$, which means $t_{lk}\in {\cal  Z}^{(0)}$.~$\Diamond$

\subsubsection {4.2 The general case $t\in {\cal  Z}^n,\  n>1$ }

Proposition \ref{p2} and Theorem~\ref{t1i} state that for $n=0,1$
the condition $t\in {\cal  Z}^n$ implies $t_{jk}\in {\cal
Z}^{n-1}$. This fact can not be extended to the case  $n>1$. If
$t\in {\cal  Z}^2$, there exist decompositions~(\ref{001}) whose
elements are not in ${\cal  Z}^1$. Nevertheless, it will be shown
that any decomposition of $t\in {\cal  Z}^n$ can be fixed to
provide  $t_{jk}\in {\cal  Z}^{n-1}$.

\begin{theo} Let $n\in\z_+$. A trigonometric polynomial $t$
belongs to ${\cal  Z}^n$  if and only if there exists a
decomposition~(\ref{001}) with $t_{jk}\in {\cal  Z}^{n-1}$,
$j,k=1\ddd d$. \label{t1} \end{theo}

{\bf Proof.} Set \ban a(x)=t({M^*}^{-1}x),
&&a_{jk}(x)=t_{jk}({M^*}^{-1}x),
\\
b_k(x)=1-\ex{{(M^*}^{-1}x,{\bf e}_k)}, &&c_k(x)=1-\ex{(x,{\bf
e}_k)}. \ean In these notation, (\ref{001}) can be rewritten as
\be b_k(x)a(x)=\sml j1d a_{jk}(x)c_j(x). \label{4} \ee Note  the
following trivial properties of $c_j$: \ba &&c_j(l)=0 \ \forall
l\in\zd, \label{5}
\\
&&D^\delta c_j\equiv0\ \forall \delta\ne r{\bf e}_j, r\in\n,
\delta\ne \nul, \label{6}
\\
&&\frac{\partial c_j}{\partial x_j}(l)\ne0\ \forall l\in\zd.
\label{7} \ea We will prove the theorem by induction on $n$. We
have a base  for $n=0$ due to Proposition~\ref{p2}. Let us prove
the inductive step: $n\to n+1$.

Assume that each polynomial $t_{jk}$ belongs to ${\cal  Z}^{n}$,
i.e. $D^\beta a_{jk}(s^*_l)=0$, $l=1\ddd m-1$, for all
$\beta\in\zd_+$, $[\beta]\le n$. Let $\alpha\in\zd$,
$[\alpha]=n+1$, $s\in D(M^*)$, $s\ne\nul$. It follows from~(\ref4)
and Leibniz formula that $$ \sum\limits_{\nul\le\beta\le\alpha}
\lll\alpha\atop\beta\rrr
 D^{\alpha-\beta}b_k(s) D^\beta a(s)=D^\alpha a_{jk}(s)c_j(s).
$$ By the inductive hypotheses, $D^\beta a(s)=0$ whenever
$\beta<\alpha$. Hence, taking into account~(\ref5), we have $$
b_k(s)D^\alpha a(s)=c_j(s)D^\alpha a_{jk}(s)=0. $$ Due to
Lemma~\ref{l4}, $b_k(s)\ne0$  for at least one $k=1\ddd d$. It
follows that $D^\alpha a(s)=0$.

Now we assume that $t\in {\cal  Z}^{n+1}$, i.e. $D^\alpha
a(s^*_\nu)=0$, $\nu=1\ddd m-1$, for all $\alpha\in\zd_+$,
$[\alpha]\le n+1$. Due to the  inductive hypotheses, there exist
polynomials $t_{jk}\in {\cal  Z}^{n-1}$ satisfying~(\ref{001}),
i.e. $D^\delta a_{jk}(s^*_\nu)=0$, $\nu=1\ddd m-1$,
 for all
$\delta\in\zd_+$, $[\delta]\le n-1$. We will construct new
trigonometric polynomials
 satisfying~(\ref{001}) and belonging to ${\cal  Z}^{n}$. Note that for $n=0$
 any trigonometric polynomials
$t_{jk}$  satisfying~(\ref{001})  belong to ${\cal  Z}^{0}$, due
to Theorem~\ref{t1i}.

Let us prove the following statement. If $D^\delta
a_{jk}(s^*_\nu)=0$, $\nu=1\ddd m-1$,
 for all $\delta\in\zd_+$, $[\delta]\le n$ and for all  $j=1\ddd l-1$ ($l=1\ddd d$),
then there exist polynomials $\tilde t_{jk}$, $j,k=1\ddd d$, such
that $\sml j1d \tilde a_{jk}c_j=\sml j1d a_{jk}c_j$, where $
\tilde a_{jk}(x)= \tilde t_{jk}({M^*}^{-1}x)$, and $D^\delta
\tilde a_{jk}(s^*_\nu)=0$, $\nu=1\ddd m-1$, $k=1\ddd d$,
 for all $\delta\in\zd_+$, $[\delta]\le n$, and for all  $j=1\ddd l$.

Let $s\in D(M^*), s\ne\nul$, $\beta\in \z_+^d$, $[\beta]=n$, $\sml
i{l+1}d\beta_i=0$. Set $\alpha=\beta+{\bf e}_l$ and note that
$[\alpha]=n+1$. It follows from~(\ref4) that \be D^\alpha\lll \sml
j1d a_{jk}(x)c_j(x)\rrr\Bigg|_{x=s}=
D^\alpha(b_k(x)a(x))|_{x=s}=0. \label8 \ee On the other hand, due
to~(\ref5), (\ref6) and the assumption of the statement,
 we have
\ban D^\alpha\lll \sml j1d a_{jk}(x)c_j(x)\rrr\Bigg|_{x=s}= \sml
j1lD^{\alpha-{\bf e}_j} a_{jk}(s)\frac{\partial c_j}{\partial
x_j}(s)= D^{\beta} a_{lk}(s)\frac{\partial c_l}{\partial x_l}(s).
\ean Combining this  with~(\ref7) and (\ref8), we get $D^\beta
a_{lk}(s) =0$.
is proved already for $l=d$ (in this case $\tilde t_{jk}=t_{jk}$,
$j,k=1\ddd d$).

Next let  $j=l+1\ddd d$, $\beta\in \z_+^d$, $[\beta]=n$, $\sml
i{l+1}d\beta_i>0$. Define the functions $$ q_{\beta
j}(x)=\frac1{-2\pi i}g_{n-1, \beta-{\bf
e}_j}(x)\sml\nu1{m-1}h_\nu(x)D^\beta a_{lk}(s^*_\nu), $$ where
$g_{N \delta}$ is a trigonometric polynomial such that
 $D^\gamma g_{N \delta}(\nul)=0$ for all $\gamma\in \z_+^d$,  $[\gamma]\le N$,
$\gamma\ne\delta$ and $D^\delta g_{N \delta}(\nul)=1$; $$
h_\nu(x)=\frac1m\sml\mu0{m-1}\ex{(x-s_\nu^*, {M}^{-1}s_\mu)}. $$
Since, by Proposition A, $h_\nu(s^*_\mu)=\delta_{\mu\nu}$, due  to
(\ref5), (\ref6) and Leibniz formula,
 it is not difficult to see that for each $\nu=1\ddd m-1$ we have
\ba &&D^\beta(c_j(x)q_{\beta j}(x))\Big|_{x=s_\nu^*}=D^\beta
a_{lk}(s^*_\nu); \label{25}
\\
&&D^\delta(c_j(x)q_{\beta j}(x))\Big|_{x=s_\nu^*}=0\ \ \ \forall
\delta\in \z_+^d,\ [\delta]\le n,\ \delta\ne\beta; \label{26}
\\
&&D^\delta(c_i(x)q_{\beta j}(x))\Big|_{x=s_\nu^*}=0\ \ \ \forall
\delta\in \z_+^d,\ [\delta]\le n-1, \forall i=1\ddd d. \label{27}
\ea Set \ban \tilde a_{lk}(x)&:=&a_{lk}(x)-\sml j{l+1}d
\sum\limits_{[\beta]=n, \
\beta_j>0\atop\beta_{l+1}=\dots=\beta_{j-1}=0} c_j(x)q_{\beta
j}(x),
\\
\tilde a_{jk}(x)&:=&a_{jk}(x)+ \sum\limits_{[\beta]=n, \
\beta_j>0\atop\beta_{l+1}=\dots=\beta_{j-1}=0} c_l(x)q_{\beta
j}(x),\ \ \ j=l+1\ddd d. \ean Because of construction, $\sml jld
\tilde a_{jk}c_j=\sml jld  a_{jk}c_j$, and, taking into
account~(\ref{25}), (\ref{26}), for each $\nu=1\ddd m-1$ we have
$D^\beta \tilde a_{lk}(s^*_\nu)=0$  whenever $[\beta]=n$, $\sml
i{l+1}d\beta_i>0$; $D^\beta \tilde a_{lk}(s^*_\nu)=D^\beta
a_{lk}(s^*_\nu)=0$  whenever $[\beta]=n$, $\sml i{l+1}d\beta_i=0$.
At last, due to~(\ref{27}), $D^\delta \tilde a_{jk}(s^*_\nu)=0$,
$j=l\ddd d$
 for all $\delta\in\zd_+$, $[\delta]\le n-1$. To complete the proof of the statement it remains
 to put $\tilde t_{jk}=\tilde a_{jk}(M^*x)$
 for $j=l\ddd d$ and $\tilde t_{jk}=t_{jk}$ for $j=1\ddd l-1$.

So, we described one step for improvement of
decomposition~(\ref{001}). Starting with $l=1$, after $(d-1)$
steps we will obtain required polynomials.~$\Diamond$

Analyzing the proof of  Theorem~\ref{t1}  it is not difficult to
describe an algorithm for finding polynomials $t_{jk}\in {\cal
Z}^{n-1}$, $j,k=1\ddd d$, in decomposition~(\ref{001}). To realize
this algorithm we will need the functions $g_{N\delta}$. Explicit
recursive formulas for
 these functions are presented in~\cite{2}.

\vspace{.5cm} \noindent {\bf ALGORITHM 2}

{\bf Input:} \ \ \ \ \  $
t\in {\cal  Z}^{N}$, $N>1$.

{\bf Output:} \ \ $
t_{jk}\in {\cal  Z}^{N-1}, j,k=1\ddd d$.

{\bf Step 1.} \ \ Using Algorithm 1, find  $t_{jk}$, $j,k=1\ddd
d$, satisfying~(\ref{001}).

{\bf Step 2.} \ \ For $n=1\ddd N-1$ do

\hspace{2.2cm} Set $t_{jk}^{(0)}:=t_{jk}$, $j,k=1\ddd d$;

\hspace{2.5cm} For each  $l=1\ddd d-1$ do

\hspace{2.8cm} For each  $k=1\ddd d$ do

$ \ \ \ \ \ a_{lk}(x)=t^{(l-1)}_{lk}({M^*}^{-1}x); $ $$ \ \ \ \
t^{(l)}_{lk}(x):=t^{(l-1)}_{lk}(x)-\sml j{l+1}d
\lll1-\ex{(M^*x,{\bf e}_j)}\rrr\sum\limits_{[\beta]=n, \
\beta_j>0\atop\beta_{l+1}=\dots=\beta_{j-1}=0}
 q_{\beta j}(M^*x);
$$

\hspace{3.1cm} For each  $j=1\ddd l$ do $$ \hspace{-8cm}
t^{(l)}_{jk}(x):=t^{(l-1)}_{jk}(x); $$

\hspace{3.1cm} For each  $j=l+1\ddd d$ do $$
t^{(l)}_{jk}(x):=t^{(l-1)}_{jk}(x)+ \lll1-\ex{(M^*x,{\bf
e}_l)}\rrr\sum\limits_{[\beta]=n, \
\beta_j>0\atop\beta_{l+1}=\dots=\beta_{j-1}=0}
 q_{\beta j}(M^*x).\Diamond
$$

In the case $d=2$, Step 2 of Algorithm 2 does not look so
frightening, it is reduced to the following.

For $n=1\ddd N-1$ do

Set  $l=1$, $t_{jk}^{(0)}:=t_{jk}$, $j,k=1,2$;

For each  $k=1,2$ do

\ban
 a_{1k}(x)&=&t^{(0)}_{1k}({M^*}^{-1}x);
\\
 t_{1k}(x)&:=&t^{(0)}_{1k}(x)-
\lll1-\ex{(M^*x,{\bf e}_2)}\rrr\sum\limits_{[\beta]=n, \
\beta_2>0}
 q_{\beta j}(M^*x);
\\
t_{2k}(x)&:=&t^{(0)}_{2k}(x)+ \lll1-\ex{(M^*x,{\bf
e}_1)}\rrr\sum\limits_{[\beta]=n, \ \beta_2>0}
 q_{\beta j}(M^*x).
\ean

In particular, for $d=2$, $N=2$, Step 2 may be realized as
follows.

Set $t_{jk}^{(0)}:=t_{jk}$, $j,k=1,2$;

 For each  $k=1,2$ do
\ban t_{1k}(x):=t^{(0)}_{lk}(x)+ \frac1{2\pi
i}\lll1-\ex{(M^*x)_2}\rrr
 \sml\nu1{m-1} h_\nu(M^*x)\frac{\partial}{\partial x_2}
t^{(0)}_{1k}({M^*}^{-1}u)\Big|_{u=s^*_\nu};
\\
t_{2k}(x):=t^{(0)}_{2k}(x)-\frac1{2\pi i} \lll1-\ex{(M^*x)_1}\rrr
\sml\nu1{m-1} h_\nu(M^*x)\frac{\partial}{\partial x_2}
t^{(0)}_{1k}({M^*}^{-1}u)\Big|_{u=s^*_\nu}. \ean

For each $n\in\n$ we introduce the set $$ \Gamma^n:=\{k\in\r^n:\ \
k_l\in\{1\ddd d\},\  l=1\ddd n\}. $$

\begin{theo}
Let  $n, n_0\in\n$, $n\le n_0$, $t\in {\cal  Z}^{n_0-1}$,
then there exist trigonometric polynomials $t_{jk}\in {\cal
Z}^{n_0-n-1}$, $k,j\in \Gamma^n$, such that \ba
\prod\limits_{l=1}^n\lll 1-\ex{(x,{\bf e}_{k_l})}\rrr t(x)&=& \sml
{j_1}1d \ldots \sml {j_n}1d t_{jk}(x)\prod\limits_{l=1}^n\lll
1-\ex{(M^*x,{\bf e}_{j_l})}\rrr, \label{13}
\\
t_{jk}(\nul)&=&\prod\limits_{l=1}^n({M}^{-1})_{j_lk_l}\,t(\nul).
\label{14} \ea \label{t2} \end{theo}

{\bf Proof.} We will proof by induction on $n$.

Base: $n=1$. Let $k\in\Gamma^1$. Due to Proposition~\ref{p2} and
Theorem~\ref{t1}, there exist trigonometric polynomials $t_{jk}\in
{\cal  Z}^{n_0-2}$, $j\in \Gamma^1$, such that~(\ref{001}) holds.
So, we have~(\ref{13})  for $n=1$, (\ref{14}) also follows from
Proposition~\ref{p2}.

Inductive step: $n-1\to n$. Let $1<n\le n_0$, $k\in\Gamma^n$,
$k^\prime:=(k_1\ddd k_{n-1})$. By the inductive hypotheses there
exist   $t_{j^\prime k^\prime}\in {\cal  Z}^{n_0-n}$, $j^\prime\in
\Gamma^{n-1}$, such that \ba \prod\limits_{l=1}^{n-1}\lll
1-\ex{(x,{\bf e}_{k_l})}\rrr t(x)= \sml {j_1}1d \ldots \sml
{j_{n-1}}1d t_{j^\prime k^\prime}(x) \prod\limits_{l=1}^{n-1}\lll
1-\ex{(M^*x,{\bf e}_{j_l})}\rrr, \label{9}
\\
t_{jk}(\nul)=\prod\limits_{l=1}^{n-1}({M}^{-1})_{j_lk_l}\,t(\nul).\hspace{4cm}
\label{10} \ea Since $n_0-n\ge0$ and the theorem is proved already
for $n=1$, for each $t_{j^\prime k^\prime}$, there exist
 trigonometric polynomials $t_{j_n k_n}\in {\cal  Z}^{n_0-n-1}$,
$j\in \Gamma^1$, such that \ba \lll 1-\ex{(x,{\bf e}_{k_n})}\rrr
t_{j^\prime k^\prime}(x)&=&
 \sml {j_n}1d t_{j_nk_n}(x)\lll 1-\ex{(M^*x,{\bf e}_{j_n})}\rrr,
\label{11}
\\
t_{j_nk_n}(\nul)&=&({M}^{-1})_{j_nk_n}t_{j^\prime k^\prime}(\nul).
\label{12} \ea Combining~(\ref9), (\ref{10}) with (\ref{11}),
(\ref{12}) we comlete the proof.~$\Diamond$

To impart a more compact form to~(\ref{13}), (\ref{14}), we
introduce the following notations. Set \ban
&&\Delta_k(x)=\Delta_k=\lll1-\ex{({\bf e}_k,x)}\rrr,\ \
\Delta(x)=\Delta=(\Delta_1\ddd \Delta_d)^T,
\\
&&\delta_k(x)=\delta_k=\lll1-\ex{({\bf e}_k,M^*x)}\rrr,\ \
\delta(x)=\delta=(\delta_1\ddd \delta_d)^T. \ean

Now Theorem~\ref{t2} can be rewritten as

\noindent {\bf Theorem 10\,$^\prime$}
{\em Let  $n, n_0\in\n$, $n\le n_0$, $t\in {\cal  Z}^{n_0-1}$,
then there exists a $d^n\times d^n$ matrix  $T$ whose entries  are
trigonometric polynomials $T_{kj}\in {\cal  Z}^{n_0-n-1}$, for
$k,j\in \Gamma^n$, such that} $$
(\Delta(x))^{[n]}t(x)=T(x)(\delta(x))^{[n]},\ \ \
T(\nul)=t(\nul)({M^*}^{-1})^{[n]}. $$

\subsection {5. Construction of masks in ${\cal Z}^n$}

A simple description of the classes ${\cal  Z}^n$ is well known in
the one-dimensional dyadic case. A general form is given by the
formula $t(x)=(1+\ex x)^nT(x)$, where $T$ is an arbitrary
trigonometric polynomial. In the multidimensional case  ${\cal
Z}^n$ can not be described in a similar way. We will give a
characterization of the class ${\cal  Z}^n$  for arbitrary
dilation matrix which allows to construct its elements in
practice.

\begin{theo} A trigonometric polynomial $t$ belongs to ${\cal
Z}^n$ if and only if the derivatives of its polyphase function
$\tau_k$, $k=0\ddd m-1$, up to order $n$ are given by
 \be
D^\alpha\tau_{k}(\nul)=\frac1{m}
\sum\limits_{\nul\le\beta\le\alpha}\lambda_\beta
\lll\alpha\atop\beta\rrr(-2\pi i r_k)^{\alpha-\beta},\ \
\alpha\in\zd_+, [\alpha]\le n, \label{23} \ee where
$\lambda_\alpha=D^\alpha t({M^*}^{-1}x)|_{x=0}$. \label{t3}
\end{theo}

{\bf Remark.} For polynomials $t$ whose polyphase functions
$\tau_0\ddd\tau_{m-1}$ form a unimodular row (i.e.  there exists a
dual row of trigonometric polynomials
$\tilde\tau_0\ddd\tilde\tau_{m-1}$ such that
$\sum_{k=0}^{m-1}\tau_k\tilde\tau_k\equiv1$), the statement of
Theorem~\ref{t3} follows from combining the results of~\cite1 and
\cite2. It was proved in these papers
 that  both conditions are equivalent to vanishing moments of
the corresponding wavelet system.

{\bf Proof.}  Assume that~(\ref{23}) holds with some  complex
numbers $\lambda_\alpha$. Let $s\in D(M^*)$.  By  Leibniz formula,
\ban D^\alpha\lll \ex{(r_k,x)}\tau_k(x)\rrr\Big|_{x=s}=
\sum\limits_{\nul\le\beta\le\alpha}\lll\alpha\atop\beta\rrr
D^\beta\lll \ex{(r_k,x)}\rrr\Big|_{x=s}
D^{\alpha-\beta}\tau_k(\nul)=
\\
\sum\limits_{\nul\le\beta\le\alpha}\lll\alpha\atop\beta\rrr\ex{(r_k,s)}
(2\pi ir_k)^\beta D^{\alpha-\beta}\tau_k(\nul)=
\\
\frac1{
m}\ex{(r_k,s)}\sum\limits_{\nul\le\beta\le\alpha}\lll\alpha\atop\beta\rrr
(2\pi ir_k)^\beta
\sum\limits_{\nul\le\gamma\le\alpha-\beta}\lambda_\gamma
\lll\alpha-\beta\atop\gamma\rrr(-2\pi i
r_k)^{\alpha-\beta-\gamma}=
\\
\frac1{ m}\ex{(r_k,s)}\sum\limits_{\nul\le\beta\le\alpha}
\sum\limits_{\nul\le\gamma\le\alpha-\beta}\lambda_\gamma
\lll\alpha-\beta\atop\gamma\rrr\lll\alpha\atop\beta\rrr(-2\pi i
r_k)^{\alpha-\gamma} \prod\limits_{j=1}^d(-1)^{-\beta_j}=
\\
\frac1{ m}\ex{(r_k,s)}\sum\limits_{\nul\le\gamma\le\alpha}
\lambda_\gamma(-2\pi i r_k)^{\alpha-\gamma}
\lll\alpha\atop\gamma\rrr
\sum\limits_{\nul\le\beta\le\alpha-\gamma}
\lll\alpha-\gamma\atop\beta\rrr
\prod\limits_{j=1}^d(-1)^{-\beta_j}. \ean Since \ban
\sum\limits_{\nul\le\beta\le\alpha-\gamma}
\lll\alpha-\gamma\atop\beta\rrr
\prod\limits_{j=1}^d(-1)^{-\beta_j}= \prod\limits_{j=1}^d
\sum\limits_{\nul\le\beta_j\le\alpha_j-\gamma_j}
\lll\alpha_j-\gamma_j\atop\beta_j\rrr(-1)^{-\beta_j}=
\\
\prod\limits_{j=1}^d(1-1)^{\alpha_j-\gamma_j}= \left\{
\begin{array}{ll} 0, & \alpha\ne\gamma,
\\
1, & \alpha=\gamma, \end{array} \right. \ean we have $$
D^\alpha\lll \ex{(r_k,x)}\tau_k(x)\rrr\Big|_{x=s}=
\frac{\lambda_\alpha}{ m}\ex{(r_k,s)}, \ \ k=0\ddd m-1. $$ It
follows from~(\ref{1}) and Proposition A that \ban
D^\alpha(t({M^*}^{-1}x)\Big|_{x=s}=
\sml k0{m-1}D^\alpha\lll \ex{(r_k,x)}\tau_k(x)\rrr\Big|_{x=s}=
\\
\frac{\lambda_\alpha}m\sml k0{m-1}\ex{(r_k,s)}= \left\{
\begin{array}{ll} \lambda_\alpha,&  \mbox{if} \ \ s=\nul,
\\
0,& \mbox{if}\ \  s\ne\nul. \end{array} \right. \ean

Now let us check that~(\ref{23}) follows from the relation $t\in
{\cal  Z}^n$. We will prove by induction on $n$. The base for
$n=0$ was established in Proposition~\ref{p1}. To prove the
inductive step $n \to n+1$, we assume that $t\in {\cal  Z}^{n+1}$
and~(\ref{23}) holds. Let $\alpha\in\zd$, $[\alpha]=n+1$, $s\in
D(M^*)$.
 By~(\ref{1}) and  Leibniz formula,
\ba D^\alpha(t({M^*}^{-1}x)\Big|_{x=s}=\sml k0{m-1}
\sum\limits_{\nul\le\beta\le\alpha}\lll\alpha\atop\beta\rrr
D^{\alpha-\beta}\lll \ex{(r_k,x)}\rrr\Big|_{x=s}
D^{\beta}\tau_k(\nul)= \nonumber
\\
\sml k0{m-1}\lll D^\alpha\tau_k(0)+
\sum\limits_{\nul\le\beta<\alpha}\lll\alpha\atop\beta\rrr
D^{\alpha-\beta}(2\pi ir_k)^{\alpha-\beta}
D^{\beta}\tau_k(\nul)\rrr\ex{(r_k,s)}. \label{24} \ea Because of
the inductive hypotheses, we have \ban
\sum\limits_{\nul\le\beta<\alpha}\lll\alpha\atop\beta\rrr (2\pi
ir_k)^{\alpha-\beta} D^{\beta}\tau_k(\nul)=
\\
\frac1m\sum\limits_{\nul\le\beta<\alpha}\sum\limits_{\nul\le\gamma\le\beta}
\lambda_\gamma\lll\alpha\atop\beta\rrr\lll\beta\atop\gamma\rrr
(2\pi ir_k)^{\alpha-\beta}(-2\pi ir_k)^{\beta-\gamma} =
\\
\frac1m\sum\limits_{\nul\le\gamma<\alpha} \lambda_\gamma (2\pi
ir_k)^{\alpha-\gamma}
\sum\limits_{\gamma\le\beta<\alpha}\lll\alpha\atop\beta\rrr\lll\beta\atop\gamma\rrr
\prod\limits_{j=1}^d(-1)^{\beta_j-\gamma_j} =
\\
\frac1m\sum\limits_{\nul\le\gamma<\alpha}
\lll\alpha\atop\gamma\rrr\lambda_\gamma (2\pi
ir_k)^{\alpha-\gamma}
\sum\limits_{\nul\le\delta<\alpha-\gamma}\lll\alpha-\gamma\atop\delta\rrr
\prod\limits_{j=1}^d(-1)^{\delta_j} =
\\
\frac1m\sum\limits_{\nul\le\gamma<\alpha}
\lll\alpha\atop\gamma\rrr\lambda_\gamma (2\pi
ir_k)^{\alpha-\gamma}
\lll\prod\limits_{j=1}^d(1-1)^{\alpha_j-\gamma_j}-\prod\limits_{j=1}^d(-1)^{\alpha_j}\rrr=
\\
-\frac1m\sum\limits_{\nul\le\gamma<\alpha}\lll\alpha\atop\gamma\rrr\lambda_\gamma
(-2\pi ir_k)^{\alpha-\gamma}. \ean Combining this with~(\ref{24}),
we obtain $$ \sml k0{m-1}\ex{(r_k,s)}\lll D^\alpha\tau_k(0)-
\frac1m\sum\limits_{\nul\le\gamma<\alpha}\lll\alpha\atop\gamma\rrr\lambda_\gamma
(-2\pi
ir_k)^{\alpha-\gamma}\rrr=D^\alpha(t({M^*}^{-1}x))\Big|_{x=s}. $$
Set $ \lambda_\alpha=D^\alpha(t({M^*}^{-1}x))\big|_{x=0}. $ Due to
Proposition A, the linear system $$ \sml
k0{m-1}\ex{(r_k,s^*_l)}y_k=\lambda_\alpha\delta_{0l},\ \ l=0\ddd
m-1, $$ has a unique solution $y_k=\frac{\lambda_\alpha}m$,
$k=0\ddd m-1$. It follows that $$ D^\alpha\tau_k(0)-
\frac1m\sum\limits_{\nul\le\gamma<\alpha}\lll\alpha\atop\gamma\rrr\lambda_\gamma
(-2\pi ir_k)^{\alpha-\gamma}=\frac{\lambda_\alpha}m, $$ which was
to be proved. $\Diamond$

So, if we want a polynomial $t$ to belong to ${\cal  Z}^n$, its
polyphase functions should have derivatives at the origin given
by~(\ref{23}). This can be easily realized for an arbitrary set of
parameters $\lambda_\beta$, $[\beta]\le n$. General forms for all
such polynomials $t$ are presented in~\cite{2}.

\subsection {6. Subdivision schemes}

In this section we apply the decompositions of Section 4 to the
analysis of convergence and smoothness of multivariate subdivision
schemes associated with general dilation matrices. We first
present basic definitions, and prove an important observation
about the matrices in a decomposition of $t\in{\cal Z}^1$.

\subsubsection {6.1  Preliminaries} Let $A_\alpha$,
$\alpha\in\zd$, be $N\times N^\prime$ matrices such that
$A_\alpha$ is  non-zero only for a finite number of $\alpha$, then
$T(x):=\sum_{\alpha\in\zd}A_\alpha \exm{(\alpha,x)}$  is a
$N\times N^\prime$ matrix whose entries are trigonometric
polynomials ($N\times N^\prime$ trigonometric matrix). The
subdivision operator $S_T=S_{T,M}$ associated with $T$ and with a
$d\times d$ dilation matrix  $M$ is defined on
$\ell_\infty^N=\ell_\infty^N(\zd)$ by $$
(S_Tf)_\alpha=\sum\limits_{\beta\in\zd}A_{\alpha-M\beta}f_\beta, \
\ f\in \ell^N_\infty. $$ It is clear that $S_T$ is a linear
bounded operator taking $\ell_\infty^N$ into
$\ell_\infty^{N^\prime}$. If $N=N^\prime$, the   operators
$S_T^n$, $n=1,2,\dots$, are well defined, and the sequence
$\{S_T^n\}_{n=1}^\infty$ is a $N\times N$ {\em matrix subdivision
scheme} ({scalar subdivision scheme if $N=1$}).

Let $t$ be a trigonometric polynomial, $T$ be a $d\times d$
trigonometric matrix,
set $L(x)=(L_1(x)\ddd L_d(x))^T$, where $L_k(x)=(1-\ex{(x,{\bf
e}_k)})t(x)$, $R(x)=(R_1(x)\ddd R_d(x))^T$, where $R_k(x)=\sml
j1d T_{kj}(x)(1-\ex{(M^*x,{\bf e}_j)})$. To each $f\in
\ell_\infty$ assign the vector-valued sequence
 $\triangledown f\in\ell_\infty^{d}$  defined by
 $(\triangledown f)_\alpha=
 (f_\alpha-f_{\alpha-{\bf e}_1}\ddd f_\alpha-f_{\alpha-{\bf e}_d})^T$.
It is clear that $S_Lf=\triangledown S_t(f)$,
$S_Rf=S_T(\triangledown f)$ for any $f\in \ell_\infty$. Hence
equality~(\ref{001}), with $t_{jk}=T_{kj}$, may be rewritten in
the form
\be
\triangledown S_t(f)=S_T(\triangledown f),\ \ \forall
f\in \ell_\infty. \label{001sd}
\ee
 Similarly, to each  $f\in \ell^N_\infty$ we assign
 $\triangledown f\in\ell_\infty^{Nd}$ defined by
$$ (\triangledown f)_\alpha= \lll(\triangledown(f,{\bf e}_1)
)_\alpha^T\ddd\triangledown(f,{\bf e}_N) )_\alpha^T \rrr^T. $$

If all  trigonometric polynomials $t_{jk}$ in the right hand side
of~(\ref{001}) are in ${\cal Z}^0$, we can decompose them (see
Theorem~\ref{t2}). This second step of the decomposition may be
rewritten as \be \triangledown S_T(f)=S_Q(\triangledown f),\ \
\forall f\in \ell^d_\infty, \label{002sd} \ee where $Q$ is a
$d^2\times d^2$ trigonometric  matrix.

\begin{prop} Let $t$ be a trigonometric polynomial, $t\in{\cal
Z}^1$, $t(\nul)=m$, and let $T$ be a $d\times d$  trigonometric
matrix satisfying~(\ref{001sd}),
$T(x):=\sum_{\alpha\in\zd}A_\alpha \exm{(\alpha,x)}$. Then $$
\sum_{\beta\in\zd}A_{\alpha-M\beta}={M^*}^{-1} \ \ \ \forall
\alpha\in\zd. $$ \label{p3} \end{prop}

{\bf Proof.} First of all, we note that
$\sum_{\beta\in\zd}A_{\alpha-M\beta}=\sum_{\beta\in\zd}A_{\alpha^\prime-M\beta}$
whenever $\alpha\equiv\alpha^\prime\pmod{M}$. So, it suffices to
check that $$
X_\nu:=\sum_{\beta\in\zd}A_{s_\nu-M\beta}={M^*}^{-1}, \ \ \
\nu=0\ddd m-1. $$ Substituting $x={M^*}^{-1}s^*_k$ into the
equality $$ T(x)=\sum\limits_{\nu=0}^{m-1}\exm{(s_\nu,x)}
\sum_{\beta\in\zd}A_{s_\nu-M\beta}\exm{(M\beta,x)}, $$ we have $
T({M^*}^{-1}s^*_k)=\sum_{\nu=0}^{m-1}\ex{(s_\nu,{M^*}^{-1}s^*_k)}X_\nu
$ It follows from Proposition~\ref{p2} and Theorem~\ref{t1i} that
$$
\sum\limits_{\nu=0}^{m-1}\exm{(s_\nu,{M^*}^{-1}s^*_k)}X_\nu=m\delta_{k0}{M^*}^{-1},
\ \ \ k=0\ddd m-1. $$ This linear system has a unique solution
$X_\nu={M^*}^{-1}$, $\nu=0\ddd m-1$, because of Proposition
A.$\Diamond$

\subsubsection {6.2 Convergence}

A scalar subdivision scheme $S_t=S_{t,M}$,
associated with a dilation matrix $M$,
is called {\em uniformly convergent}, if for any $f\in
\ell_\infty$, there exists  a continuous function $S^\infty_tf$
such that \be \lim\limits_{k\to
\infty}\|S_t^kf-S^\infty_tf(M^{-k}\cdot)\|_\infty=0 \label{31} \ee
and if for at least one $f\in \ell_\infty$, the limit function
$S_t^\infty f$ is not identically zero.

The issue of the convergence of a multivariate scalar
subdivision scheme associated with a dilation matrix $M$ is
studied in ~\cite{10}. The convergence result there is limited to
dilation matrices having a self-similar tile.
A $d\times d$ dilation matrix $M$ is said to have a {\em
self-similar tile} if there exist a set of digits $D(M)$ and a
bounded set $E\subset\r^d$ whose integer translates form a
disjoint decomposition of $\rd$ such that \be
ME=\cupm\limits_{s\in D(M)}(E+s),
\label{tiling} \ee

\begin{theo}\cite{10} Let $S_t=S_{t,M}$ be a scalar subdivision
scheme and let there exists a self-similar tile related to $M$.
Then $S_t$ is uniformly convergent
if and only if: a) $t(0)=m$;
b) \ there exists a $d\times d$
 trigonometric matrix $T$ such that (\ref{001sd}) holds,
$c)\ \ \lim\limits_{k\to\infty}\sup\limits_{f\in\ell_\infty\atop
\|\triangledown f\|_\infty=1} \|S^k_T\triangledown
f\|_{\ell_\infty^d}=0.$
\label{t4}
\end{theo}

Here we present  sufficient conditions for convergence of a
scalar multivariate
subdivision scheme associated with any
dilation matrix $M$, and then give an example
of a subdivision scheme satisfying these conditions, where the
verefication of the conditions is done with Algorithm 1.

\smallskip

{\bf Theorem  $\bf 13^\prime$} {\em
Let $t\in {\cal Z}^0$ and  $M$   be a $d\times d$ dilation
matrix. Then there exists a trigonometric matrix $T$ of order
$d\times d$  satisfying
(\ref{001sd}).  Moreover, if   $t(0)=m$ and
\be
\lim_{k\to+\infty}\|S_T^k\|_\infty=0,
\label{601}
\ee
then the subdivision scheme  $S_t$  is uniformly convergent.}

\smallskip

First we prove a simple lemma.

\begin{lem}
\label{lastlemma1}
Let $n\in\z_+$, $M$   be a $d\times d$ dilation
matrix, and let  $\phi$ be a compactly supported
function   satisfying the refinement equation
\be
\phi(x)=\sum_{\alpha\in \zd} a_\alpha\phi(Mx-\alpha),\ \ \forall
x\in\r^d.
\label{refinement1}
\ee
Then
 \be
 \sum_{\beta\in\zd}f_\beta\phi(x-\beta)=
 \sum_{\alpha\in\zd}(S_{\tilde t I_n}~f)_\alpha\phi(Mx-\alpha),\ \
 \forall x\in\r^d \label{adjoint}
 \ee
 for any $f\in\ell^n_\infty$, where
 $${\tilde t}(x)=\sum_{\alpha\in\zd}a_\alpha \exm{(x,\alpha)}$$
\end{lem}

{\bf Proof}.
Using~(\ref{refinement1}), we have
\ban\sum_{\alpha\in\zd}(S_{\tilde t I_n}f)_\alpha\phi(Mx-\alpha)=
\sum_{\alpha\in\zd}\sum_{\beta\in\zd}a_{\alpha-M\beta}f_\beta\,\phi(Mx-\alpha)=
\\
\sum_{\beta\in\zd}f_\beta
\sum_{\alpha\in\zd}a_\alpha\phi(M(x-\beta)-\alpha)=\sum_{\beta\in\zd}f_\beta
\phi(x-\beta).\Diamond
\ean

\smallskip

 {\bf Proof of Theorem $13^\prime$}. The existence of $T$
satisfying (\ref{001sd})  follows from
Proposition~\ref{p2}, as indicated in Subsection 6.1.

Let $\phi$ be a continuous compactly supported function
satisfying (\ref{refinement1}) and the interpolatory conditions
\be
\phi(\alpha)=\delta_{\alpha\nul},\ \ \forall \alpha\in\zd,
\label{partition1}
\ee
The existence of such a function for an
arbitrary  dilation matrix $M$ is proved in~\cite[Proposition 4.1]{15}.
It is proved in \cite{15} that the corresponding
mask $ {\tilde t}$
is in ${\cal Z}^0$ and that $\tilde t(0)=m$.
To prove the second part of the claim,
we will show that the sequence of  functions
$$
F_k(x)=\sum_{\alpha\in\zd}\phi(M^kx-\alpha)(S^k_tf)_\alpha,\ \ \ k=1,2,\dots,
$$
is a Cauchy sequence in
$C(\rd)$ for any $f\in\ell_\infty$.

By Lemma~\ref{lastlemma1} with $n=1$, we get for $f\in\ell_\infty$
$$
F_k(x)=\sum_{\alpha\in\zd}\phi(M^{k}x-\alpha)
(S^k_tf)_\alpha=\sum_{\alpha\in\zd}\phi(M^{k+1}x-\alpha)
(S_{\tilde t }S^k_t f)_\alpha,
$$
and therefore
\ba
F_{k+1}(x)-F_k(x)=
\sum_{\alpha\in\zd}\phi(M^{k+1}x-\alpha)
((S_t-S_{\tilde t })S^k_tf)_\alpha.
 \label{631}
 \ea

Since $t$ and $\tilde t$ are in
${\cal Z}^0$ and $t(0)=\tilde t(0)$,
we get from Lemma~\ref{l2} that
$$
t(x)-\tilde t(x)=\sml k1d q_k(x)\lll
1-\ex{(M^*x,{\bf e}_k)}\rrr.
$$
This leads, by similar derivations to those
leading to (25), to
the existence of a vector of subdivision schemes  $(S_{q_1}\ddd S_{q_d})$,
such that for any $f\in\ell_\infty$,
$$
(S_t-S_{\tilde t})f=S_{t-\tilde t}f=(S_{q_1}\ddd S_{q_d})\triangledown f.
$$

Hence, (\ref{631}) can be rewritten as
\ban
F_{k+1}(x)-F_k(x)=\sum_{\alpha\in\zd}
((S_{q_1}\ddd S_{q_d})\triangledown S^k_t f)_\alpha\phi(M^{k+1}x-\alpha),
\ean
Using relation~(\ref{001sd}), we obtain
\ban
|F_{k+1}(x)-F_k(x)|=\left|\sum_{\alpha\in\zd}
((S_{q_1}\ddd S_{q_d}) S^k_T\triangledown f)_\alpha\phi(M^{k+1}x-\alpha)\right|\le
\\
C\| S^k_T\triangledown f\|_\infty
\le C\| S^k_T\|_\infty|\triangledown f|,
\ean
where $C$ depends on $\phi, q_1,\ldots,q_d$.

Condition~(\ref{601}) is equivalent to the existence of a positive integer $L$
such that $\| S^L_T\|_\infty=\mu<1$.
Thus
\ban
|F_{k+1}(x)-F_k(x)|_{}\le C|\triangledown f|\mu^{\frac{k}{L}}.
\ean
This yields that for all $k,n\in\n$
\be
|F_{k+n}(x)-F_k(x)|_{}\le C |\triangledown f|\sum_{j=0}^{n-1}\mu^{\frac{k+j}{L}}\le
 \frac{C}{1-\mu} |\triangledown f|\,\mu^{\frac{k}{L}},
 \label{651}
 \ee
 which implies that $\{F_k\}$ is a Cauchy sequence. Denote the
 limit function by $F$ which  is, evidently, a continuous function.
 Passing to the limit in~(\ref{651}) as $n\to\infty$, we have
\be
\lim_{k\to\infty}\sup_{x\in\rd}|F(x)-F_k(x)|_{}=0,
\label{671}
\ee

Due to~(\ref{671}),
\ban
\lim_{k\to\infty}\sup_{x\in\rd}\left|F(x)-
\sum_{\alpha\in\zd}(S_t^k f)_\alpha
\phi(M^{k}x-\alpha)\right|_{}=0,
\ean
Substituting
$x=M^{-k}\beta$, $\beta\in\zd$, we have
$$
\lim_{k\to\infty}\sup_{\beta\in\zd}\left|F(M^{-k}\beta)-
\sum_{\alpha\in\zd}(S_t^k f)_\alpha\phi(\beta-\alpha)\right|_{}=0,
$$
The claim of the
theorem follows now from (\ref{partition1}). ~$\Diamond$

\smallskip

{\bf Example.\ \ \ }
Let
$
M=\left(\begin{array}{rr} 0 & 2 \\ 2 & -1
\end{array} \right),
$
$$
s_0=\left(\begin{array}{r} 0 \\ 0\end{array} \right),
s_1=\left(\begin{array}{r} 1 \\ 0\end{array} \right),
s_2=\left(\begin{array}{r} 0 \\ 1\end{array} \right),
s_3=\left(\begin{array}{r} 1\\ 1 \end{array} \right).
$$
Define the polyphase functions of  a polynomial $t$ by

\ban
\tau_0(x)&=&\frac1{16}(4+4z_1+4z_2+4z_1z_2),
\\
\tau_1(x)&=&\frac1{16}(5+4z_1+z_1^{-1}+2z_2+3z_2^{-1}+z_1z_2),
\\
\tau_2(x)&=&\frac1{16}(4+z_1+2z_1^{-1}+5z_2+z_2^{-1}+3z_1z_2),
\\
\tau_3(x)&=&\frac1{16}(5+z_1+4z_1^{-1}+z_2+3z_2^{-1}+z_1z_2+z_1^{-1}z_2).
\ean
It is not difficult to see that $t\in {\cal  Z}^0$,
and that $t(0)=m$.
Using Algorithm 1, we find
$$
\begin{array}{ll}
\tau_{110}(x)=\frac1{16}(-1+2z_2+z_1z_2+2z_2^2+z_1z_2^2),
&
\tau_{210}(x)=\frac1{16}(5+3z_2),
\\
\\
\tau_{111}(x)=\frac1{16}(z_1^{-1}+3z_2),
&
\tau_{211}(x)=\frac1{16}(5+3z_2^{-1}),
\\
\\
\tau_{112}(x)=\frac1{16}(-1+2z_1^{-1}+z_2+z_2^2-z_1^{-1}z_2),
&
\tau_{212}(x)=\frac1{16}(5+3z_2+z_2^{-1}),
\\
\\
\tau_{113}(x)=\frac1{16}(2z_1^{-1}+2z_2+z_1^{-1}z_2),
&
\tau_{213}(x)=\frac1{16}(5+2z_2^{-1}),
\\
\\
\tau_{120}(x)=\frac1{16}(1+z_1+4z_2+z_2^{-1}+3z_1z_2),
&
\tau_{220}(x)=\frac1{16}z_2^{-1},
\\
\\
\tau_{121}(x)=\frac1{16}(2+z_1+z_1^{-1}+z_2+3z_2^{-1}+z_1z_2),
&
\tau_{221}(x)=0,
\\
\\
\tau_{122}(x)=\frac1{16}(3+z_2+2z_1^{-1}),
&
\tau_{222}(x)=\frac1{16}z_2^{-1},
\\
\\
\tau_{123}(x)=\frac1{16}(3+3z_1^{-1}+z_1^{-1}z_2),
&
\tau_{223}(x)=0.
\end{array}
$$
By Theorem $13^\prime$, to prove that the subdivision scheme $S_t$
is convergent,
it remains to check that the norm of the matrix subdivision
operator $S_T$  is strictly less than 1.
This requirement is fulfilled because
$$
\|S_T\|_\infty\le\max\limits_{(\nu,k)}
(\|\tau_{1k\nu}\|_{\ell_1}+\|\tau_{2k\nu}\|_{\ell_1})\le\frac{15}{16}.
$$
(Here we identify a trigonometric polynomial with its
sequence of Fourier coefficients).

\subsubsection{6.3 Smoothness}

Now we discuss how to study the smoothness of the limit function
of a uniformly convergent scalar subdivision scheme.

Let $t=\sum\limits_{\alpha\in\zd}a_\alpha \exm{(\alpha,\cdot)}$ be
a trigonometric polynomial, and let  $\tau_{\nu}$, $\nu=0\ddd
m-1$, be its polyphase functions. It is well known (see,
e.g.,~\cite{10})  that $$
\tau_{\nu}(0)=\sum_{\beta\in\zd}a_{s_\nu-M\beta}=1,\  \nu=0\ddd
m-1, $$ whenever~(\ref{31}) is fulfilled for at least one
$f\in\ell_\infty$ for which $S_t^{\infty}f \not\equiv 0$.
So, if $S_t$ uniformly converges  then, due to
Propositions~\ref{p1} and \ref{p2}, there exists a $d\times d$
trigonometric matrix $T$  such that~(\ref{001sd}) holds.

 Following ~\cite{14}, we introduce the notions of a {\em
normalized subdivision scheme} and its {\em  subconvergence}.

For  $X$  a $N\times N$   matrix, and  $T$  a $N\times N$
trigonometric  matrix, the sequence
$\{{X^*}^kS_T^k\}_{k=1}^\infty$ is called the {\em normalized ( by
$X$) subdivision scheme} $S_T$.

We say that a subdivision scheme $S_T=S_{T,M}$  normalized by $X$
is {\em uniformly convergent} on a subspace $L$ of
$\ell^N_\infty$, $L\ne \{0\}$, if for any $F\in L$, there exists a
continuous vector-valued function $S^\infty_{T/X}F$
($S^\infty_{T/X}F(x)\in\r^N$) such that \be \lim\limits_{k\to
\infty}\|{X^*}^kS_T^kF-S^\infty_{T/X}F(M^{-k}\cdot)\|_\infty=0.
\label{33} \ee

\noindent and if for at least one $F\in L$ the limit $S_{T/X}F$ is
not identically zero.

We say that a subdivision scheme $S_T$  normalized by $X$ is {\em
uniformly subconvergent on a subspace $L$},  $L\subseteq
\ell^d_\infty$, $L\ne \{0\}$, if for some infinite set $\frak
N\subset \n$ and for any $F\in L$ there exists  a continuous
vector-valued function $S^\infty_{T/X}F$
($S^\infty_{T/X}F(x)\in\r^d$) such that \be \lim\limits_{k\in\frak
N\atop k\to
\infty}\|{M^*}^kS_T^kF-S^\infty_{T/X}F(M^{-k}\cdot)\|_\infty=0.
\label{33} \ee and if for at least one $F\in L$, the limit
$S_{T/X}F$ is not identically zero.

Next, following~\cite{13}, we  consider the class of isotropic
dilation matrices. A matrix $M$ is called {\em isotropic} if there
exists a constant $C$ such that
$\|M^k\|_\infty\|M^{-k}\|_\infty\le C$ for all $k\in\n$.

 Here we give a detailed proof of the sufficiency of a
necessary and sufficirnt condition for the limits of a scalar
subdivision scheme to be in $C^1$. This condition is stated in
~\cite{14} with a sketch of a proof, which is not applicable for
general isotropic dilation matrices. Indeed, by Proposition
\ref{p3}  here, $S_T$ cannot be of full rank, as required by the
sketch of the proof, whenever  $M$ is not a multiple of the
identity matrix.

\begin{theo} Let $M$ be isotropic, let a scalar subdivision scheme
$S_t=S_{t,M}$ be uniformly convergent,
and let $S_T$ be a   $d\times d$ matrix subdivision scheme
satisfying~(\ref{001sd}). If the subdivision scheme $S_T$
normalized by $M$  uniformly subconverges on    $\triangledown
\ell_\infty$, then  the  function $S^\infty_tf$ is in $C^1(\rd)$
for all  $f\in \ell_\infty$. \label{t6} \end{theo}

 First we prove a lemma.
 \begin{lem} Let  $\frak M$ be an
infinite subset of $\n$, $g_k\in\ell_\infty^n$, $k\in \frak M$,
$g$ be a continuous vector-valued function ($g(x)\in\r^n$), and
let $\Psi_k$, $k\in \frak M$,  be  $n\times n$ matrix-valued
 functions which are uniformly bounded, uniformly compactly supported and such that
\be \sum_{\alpha\in\z^n}\Psi_k(\cdot-\alpha)\equiv I_n, \ \ \
\forall k\in \frak M. \label{38} \ee
 If
\be \lim_{k\in\frak M\atop k\to
\infty}\|g_k-g(M^{-k}\cdot)\|_\infty=0, \label{35} \ee then the
sequence of  vector-valued functions
$\sum_{\alpha\in\z^d}\Psi_k(M^k\cdot-\alpha)(g_k)_\alpha$, $k\in
\frak M$, uniformly converges to g on any compact set
$K\subset\r^d$. \label{l6} \end{lem}

{\bf Proof}. Let $K\in\rd$ be a compact set, $x\in K$, $k\in\frak
M$. Because of~(\ref{38}), \ba
\sum_{\alpha\in\z^n}\Psi_k(M^kx-\alpha)(g_k)_\alpha-g(x)=
\sum_{\alpha\in\z^n}\Psi_k(M^kx-\alpha)\lll(g_k)_\alpha-g(x)\rrr=
\nonumber
\\
\sum_{\alpha\in\Omega(M^kx)}\Psi_k(M^kx-\alpha)\lll(g_k)_\alpha-g(M^{-k}\alpha)\rrr+
\hspace{3cm} \nonumber
\\
\sum_{\alpha\in\Omega(M^kx)}\Psi_k(M^kx-\alpha)\lll
g(M^{-k}\alpha)-g(x)\rrr, \label{41} \ea
where $\Omega(t)=\cupm\limits_{k\in\frak M}\Omega_k(t)$, with $
\Omega_k(t)=\{\alpha\in\zd:\ \ \Psi_k(t-\alpha)\ne{\Bbb O}_n\}$.

It is clear, that $\sharp\,\Omega(t)\le C$, where $C$ is a
constant depending only on the joint support of the functions
$\Psi_k$. It follows from~(\ref{35}) that \be \lim_{k\in\frak
M\atop k\to \infty}\sup_{x\in K}
\left|\sum_{\alpha\in\Omega(M^kx)}\Psi_k(M^kx-\alpha)
\lll(g_k)_\alpha-g(M^{-k}\alpha)\rrr\right|_{}=0. \label{39} \ee
If  $\Psi_k(M^kx-\alpha)\ne{\Bbb O}_n$, then
$|M^{-k}\alpha-x|\le\|M^{-k}\|R$, where $R$ is the radius of a
ball containing the supports of all $\Psi_k$. Due to the uniform
continuity of $g$ on a compact set, this yields \be
\lim_{k\in\frak M\atop k\to \infty}\sup_{x\in K}
\left|\sum_{\alpha\in\Omega(M^kx)}\Psi_k(M^kx-\alpha) \lll
g(M^{-k}\alpha)-g(x)\rrr\right|_{}=0. \label{40} \ee To complete
the proof it remains to combine~(\ref{39}) and (\ref{40}) with
(\ref{41}).$\Diamond$

{\bf Proof of Theorem~\ref{t6}}.  For $\sigma\in\n^d$, the
function $$b_\sigma(x)
=\prod\limits_{j=1}^d(\chi_{[0,1]}*\dots*\chi_{[0,1]})(x_j)$$ is
the tensor product B-spline of order $\sigma$, when the number of
convolutions in the above product are $\sigma_j+1$. In the
following we take $\sigma_j>1,\ j=1\ddd d$. Thus, $b_\sigma$ has
the following properties: \be
\sum_{\alpha\in\zd}b_\sigma({x-\alpha})\equiv1, \label{36} \ee \ba
\mbox{grad}\, b_\sigma=
\Big( b_{\sigma-{\bf e}_1}(\cdot-{\bf e}_1)-b_{\sigma-{\bf
e}_1}\ddd b_{\sigma-{\bf e}_d}(\cdot-{\bf e}_d)-b_{\sigma-{\bf
e}_d}\Big)^T. \label{37} \ea Given $f\in L$, using~(\ref{36}) and
Lemma~\ref{l6} with $\frak M=\n$, $n=1$, $g_k=S_t^kf$,
$g=S^\infty_tf$, we obtain that the sequence
$F^k:=\sum_{\alpha\in\z^n}b_\sigma(M^k\cdot-\alpha)(S^k_tf)_\alpha$
uniformly converges to $S^\infty_tf$ on any compact set
$K\subset\rd$.  These functions  are in  $C^1(\rd)$ by the choice
of $\sigma$. To prove the theorem  it remains to check that the
sequence of vector-valued functions $\mbox{grad}\, F^k$ uniformly
converges to a continuous vector-valued function on any compact
set $K\subset\rd$.

Set
 $$ B_\sigma=\left(\begin{array}{rrrr} b_{\sigma-{\bf e}_1}
&\dots&0
\\
\vdots&\ddots&\vdots
\\
0&\dots &b_{\sigma-{\bf e}_d} \end{array} \right)
$$ Using~(\ref{001sd}), we have \ban \mbox{grad}\,
\sum_{\alpha\in\zd}b_\sigma(M^kx-\alpha)(S^k_tf)_\alpha=
{M^*}^k\sum_{\alpha\in\zd}\mbox{grad}\,
b_\sigma(M^kx-\alpha)(S^k_tf)_\alpha=
\\
{M^*}^k\sum_{\alpha\in\zd}B_\sigma(M^kx-\alpha)(\triangledown
S^k_tf)_\alpha=
\sum_{\alpha\in\zd}{M^*}^kB_\sigma(M^kx-\alpha){M^*}^{-k}({M^*}^kS^k_T\triangledown
f)_\alpha. \ean Due to the normalized uniform subconvergence of
$S_T$ on $\triangledown L$, there exists an  infinite set $\frak
N\subset \n$ such that for any $F\in \triangledown L$ there exists
a continuous vector-valued function $S^\infty_{T/M}F$  for
which~(\ref{33}) holds. Let $f\in L$. Since
$\sum_{\alpha\in\zd}{M^*}^kB_\sigma(x-\alpha){M^*}^{-k}\equiv
I_d$, by Lemma~\ref{l6} with $\frak M=\frak N$, $n=d$,
$g_k={M^*}^kS^k_T\triangledown f$, $g=S^\infty_{T/M}\triangledown
f$ and $\Psi_k={M^*}^kB_\sigma{M^*}^{-k}$ (the functions $\Psi_k$
are uniformly bounded because $M$ is isotropic),  we can state
that the sequence $\mbox{grad}\, F^k$, $k\in \frak N$,  uniformly
converges to $S^\infty_{T/M}\triangledown f$ on any compact set
$K\subset\rd$.~$\Diamond$

 A direct conclusion from Theorem~\ref{t6} in case equation
~(\ref{001sd}) holds for a trigonometric polynomial $t$, is that
the limit function of the subdivision scheme $S_t$ is in $C^1$ if
$S_t$ is uniformly convergent and if $S_T$ normalized by $M$ is
uniformly subconvergent on $\triangledown \ell_\infty$.

 How to check the uniform convergence of a scalar subdivision
scheme is discussed in Subsection~6.2. The method is based on
Theorem~\ref{t4} from \cite{10}.
In the following we prove a sufficient condition for the uniform
convergence of $S_T$ normalized by $M$ on $\triangledown
\ell_\infty$,
%
because to the best of our knowledge, such a theorem has not been
published yet.
Here we formulate and prove such a theorem  for our specific case,
namely for checking the smoothness of limits of a convergent
scalar multivariate subdivision scheme, corresponding to a
trigonometric polynomial in ${\cal Z}^1$.

\begin{theo} Let $t\in {\cal Z}^1$, $M$   be a $d\times d$ dilation
matrix. Then there exist trigonometric matrices $T,Q$ of orders
$d\times d$ and $d^2\times d^2$ respectively, satisfying
(\ref{001sd}) and (\ref{002sd}). Moreover, if  $t(0)=m$ and
\be
\|{M^*}^L\|_\infty\|S_Q^L\|_\infty\le\mu<1 \label{331} \ee for
some positive integer  $L$,
then the subdivision scheme  $S_T$ normalized by $M$  uniformly
converges on $\triangledown\ell_\infty$.
\label{lastheorem}
\end{theo}

{\bf Proof}. The existence of $T,Q$,
satisfying (\ref{001sd}) and  (\ref{002sd}), follows from
Theorem~\ref{t2}, as indicated in Subsection 6.1.

To prove the second part of the claim,
we will show that the sequence of vector-valued functions
$$
F_k(x)=\sum_{\alpha\in\zd}\phi(M^kx-\alpha)({M^*}^{k}S^k_T\triangledown f)
_\alpha,\ \ \
k=1,2,\dots,
$$
with $\phi$ as in the proof of Theorem $13^\prime$,
is a Cauchy sequence in
$(C(\rd))^d$ for any $f\in\ell_\infty$.

Recall that $\phi$ is a continuous compactly supported function
satisfying (\ref{refinement1}) and the interpolatory conditions
(\ref{partition1}), and that
${\tilde t}$, defined  as in Lemma \ref{lastlemma1},
is in ${\cal Z}^0$ and satisfies $\tilde t(0)=m$.


By Lemma~\ref{lastlemma1} with $n=d$,
$$
F_k(x)={M^*}^{k}\sum_{\alpha\in\zd}\phi(M^{k}x-\alpha)
(S^k_T\triangledown
f)_\alpha={M^*}^{k}\sum_{\alpha\in\zd}\phi(M^{k+1}x-\alpha)
(S_{\tilde t I_d}S^k_T\triangledown f)_\alpha.
$$
This yields that
\ba
F_{k+1}(x)-F_k(x)=
{M^*}^{k}\sum_{\alpha\in\zd}\phi(M^{k+1}x-\alpha)
(({M^*}S_T-S_{\tilde t I_d})S^k_T\triangledown f)_\alpha=
\nonumber
\\
{M^*}^{k}\sum_{\alpha\in\zd}\phi(M^{k+1}x-\alpha) (S_{\widetilde
T}S^k_T\triangledown f)_\alpha, \label{431}
\ea
 where $\widetilde
T={M^*}T-\tilde t I_d$. The entries of the matrix $\widetilde T$
are trigonometric polynomials in ${\cal Z}^0$ because $\tilde
t\in{\cal Z}^0$ as was mentioned above, and  the  entries of $T$
are in ${\cal Z}^0$ due to Theorem~\ref{t2}. Moreover, by
Proposition~\ref{p3},
$$
\widetilde T(\nul)={M^*}T(\nul)-\tilde
t(\nul) I_d= {M^*}m{M^*}^{-1}-m I_d={\Bbb O}_d.
$$
Due to
Lemma~\ref{l2}, it follows that for any $g\in \ell^d_\infty$ \be
S_{\widetilde T}g=S_{\widetilde Q}\triangledown g, \label{421} \ee
where $\widetilde Q$ is a $d\times d^2$  trigonometric matrix.

So, (\ref{431}) may be rewritten as \ban
F_{k+1}(x)-F_k(x)={M^*}^{k}\sum_{\alpha\in\zd}\phi(M^{k+1}x-\alpha)
(S_{\widetilde Q}\triangledown S^k_T\triangledown f)_\alpha, \ean
and \be |F_{k+1}(x)-F_k(x)|_{}\le C\|{M^*}^{k}\|_\infty
\|S_{\widetilde Q}\triangledown S^k_T\triangledown f\|_\infty,
\label{441} \ee where $C$ depends only on $\phi$.

It follows  from~(\ref{002sd}) that $\triangledown S^k_T g= S^k_Q
\triangledown g$ for any $g\in \ell_\infty^d$, and hence \be
\triangledown S^k_T \triangledown f= S^k_Q \triangledown^2 f.
\label{511} \ee

Thus we get from (\ref{441}), in view of (\ref{331}),

\ban |F_{k+1}(x)-F_k(x)|_{}\le C\|{{
M}^*}^{k}\|_\infty\|S_{\widetilde Q}\|_\infty \|S^k_Q
\triangledown^2 f\|_\infty\le\hspace{2cm} \nonumber
\\
C_1 {\|{{ M}^*}^{L}\|_\infty}^{\frac{k}{L}}\, {\|S^L_Q\|_\infty
}^{\frac{k}{L}}\,|\triangledown^2 f|\le C_1\,|\triangledown^2 f|
\,\mu^{\frac{k}{L}}.
\ean
This yields that for all $k,n\in\n$
\be
|F_{k+n}(x)-F_k(x)|_{}\le C_1\,\|\triangledown^2
f\|_\infty \sum_{j=0}^{n-1}\mu^{\frac{k+j}{L}}
\le C_2\,\|\triangledown^2
f\|_\infty\,\mu^{\frac{k}{L}},
\label{451}
\ee
which implies that $\{F_k\}$ is a Cauchy sequence. Denote the
limit vector-valued function by $F$ which  is, evidently continuous.
Passing to the limit in~(\ref{451}) as $n\to\infty$, we have \be
\lim_{k\to\infty}\sup_{x\in\rd}|F(x)-F_k(x)|_{}=0, \label{471} \ee

Due to~(\ref{471}),
\ban
\lim_{k\to\infty}\sup_{x\in\rd}\left|F(x)-
\sum_{\alpha\in\zd}({M^*}^kS_T^k\triangledown f)_\alpha
\phi(M^{k}x-\alpha)\right|_{}=0,
\ean
Substituting
$x=M^{-k}\beta$, $\beta\in\zd$, we have
$$
\lim_{k\to\infty}\sup_{\beta\in\zd}\left|F(M^{-k}\beta)-
\sum_{\alpha\in\zd}({M^*}^kS_T^k\triangledown
f)_\alpha\phi(\beta-\alpha)\right|_{}=0,
$$
The claim of the
theorem follows now from (\ref{partition1}). ~$\Diamond$

 \end{document}